\documentclass[a4paper,12pt]{amsart}

\usepackage{float}
\usepackage{euscript,eufrak,verbatim}
\usepackage{graphicx}
\usepackage[usenames]{color}
\usepackage[colorlinks,linkcolor=red,anchorcolor=blue,citecolor=blue]{hyperref}
\usepackage{amsmath}
\usepackage{amsthm}
\usepackage[all]{xy}
\usepackage{graphicx} 
\usepackage{amssymb}

\theoremstyle{plain}
\newtheorem{main theorem}{Main Theorem}
\newtheorem{theorem}{Theorem}[section]
\newtheorem{lemma}[theorem]{Lemma}

\newtheorem{corollary}[theorem]{Corollary}
\newtheorem{proposition}[theorem]{Proposition}

\theoremstyle{definition}

\newtheorem{remark}[theorem]{Remark}
\newtheorem{example}[theorem]{Example}

\newtheorem{problem}[theorem]{Problem}

\makeatletter 
\@addtoreset{equation}{section}
\numberwithin{equation}{section}

\newcommand{\norm}[1]{\left\lVert#1\right\rVert}

\newcommand{\diam}{\mathrm{Diam}}

\newcommand{\mdim}{\mathrm{mdim}}
\newcommand{\vol}{\mathrm{vol}}
\newcommand{\widim}{\mathrm{Widim}}
\newcommand{\htop}{h_{\mathrm{top}}}
\newcommand{\umdimm}{\overline{\mathrm{mdim}}_{\mathrm{M}}}
\newcommand{\lmdimm}{\underline{\mathrm{mdim}}_{\mathrm{M}}}
\newcommand{\mdimm}{\mathrm{mdim}_{\mathrm{M}}}
\newcommand{\ldimm}{\underline{\dim}_{\mathrm{M}}}
\newcommand{\udimm}{\overline{\dim}_{\mathrm{M}}}

\title{Application of waist inequality to entropy and mean dimension}

\author{Ruxi Shi}

\address{Sorbonne Universite, LPSM, 75005 Paris, France}

\email{ruxi.shi@upmc.fr}

\author{Masaki Tsukamoto}

\address
{Department of Mathematics, Kyoto University, Kitashirakawa Oiwake-cho, Sakyo-ku, Kyoto 606-8502, Japan}

\email{tsukamoto@math.kyoto-u.ac.jp}

\begin{document}



\thanks{R.S was partially supported by Fondation Sciences Mathématiques de Paris. M.T. was supported by JSPS KAKENHI JP21K03227.}

\maketitle

\begin{abstract}
Waist inequality is a fundamental inequality in geometry and topology.
We apply it to the study of entropy and mean dimension of dynamical systems.
We consider equivariant continuous maps $\pi: (X, T) \to (Y, S)$ between dynamical systems and 
assume that the mean dimension of the domain $(X, T)$ is larger than the mean dimension of the target $(Y, S)$.
We exhibit several situations for which the maps $\pi$ necessarily have positive conditional metric mean dimension.
This study has interesting consequences to the theory of topological conditional entropy.
In particular it sheds new light on a celebrated result of Lindenstrauss and Weiss about minimal dynamical systems
non-embeddable in $[0,1]^{\mathbb{Z}}$.
\end{abstract}

\section{Introduction} \label{section: introduction}

\subsection{Topological conditional entropy}  \label{subsection: topological conditional entropy}

Waist inequality is a deep geometric-topological inequality which incorporates the ideas of 
isoperimetric inequality with the Borsuk--Ulam theorem in a highly nontrivial way.
It states that, for $n\geq m$, every continuous map $f:S^n\to \mathbb{R}^m$ has some fiber of “large volume” in an appropriate sense.
We will review it in Section \ref{section: waist inequality} below.
Its origin goes back to the work of Almgren \cite{Almgren} on geometric measure theory.
Gromov mentioned a version of the waist inequality in his famous ``filling’’ paper \cite[Appendix I (F)]{Gromov83}.
Later he developed its full theory in \cite{Gromov03}.

Variants of the waist inequality have appeared in several interesting subjects:
convex geometry \cite{Vershynin, Giannopoulos--Milman--Tsolomitis, Klartag17},
generalizations of Morse theory \cite{Gromov09}, combinatorics \cite{Gromov10},
knots \cite{Pardon, Gromov--Guth}  and
complex analysis \cite{Klartag18}.
The purpose of this paper is to further add a new item to this list;
we explain applications of the waist inequality to \textit{topological dynamics}.

Our starting point is the theory of topological conditional entropy.
We first quickly prepare basic definitions.
Details about topological conditional entropy can be found in 
the book of Downarowicz \cite[Section 6]{Downarowicz}.

 A pair $(X, T)$ is called a \textbf{dynamical system} if $X$ is a compact metrizable space and
 $T:X\to X$ is a homeomorphism.
 Let $d$ be a metric on $X$. 
 For a natural number $N$ we define a new metric $d_N$ on $X$ by 
 \[  d_N(x, y) = \max_{0\leq n < N} d\left(T^n x, T^n y\right). \]
 For a positive number $\varepsilon$ and a subset $E\subset X$ we denote by $\#\left(E, d_N, \varepsilon\right)$ the 
 minimum number of $n$ such that there exists an open covering $\{U_1, \dots, U_n\}$ of $E$ 
 satisfying $\diam(U_i, d_N) < \varepsilon$ for all $1\leq i \leq n$.
 Here $\diam(U_i, d_N)$ is the diameter of $U_i$ with respect to the metric $d_N$.
 When $E = \emptyset$, we set $\#\left(E, d_N, \varepsilon\right) = 0$.
 
We define the \textbf{topological entropy} of $(X, T)$ by 
\[  \htop(X, T) = \lim_{\varepsilon\to 0} \left(\lim_{N\to \infty} \frac{\log \#\left(X, d_N, \varepsilon\right)}{N}\right). \]
This is a topological invariant of dynamical systems. Namely its value is independent of the choice of the metric $d$.
For a closed subset $E\subset X$ we also define
\[  \htop(E, T) = \lim_{\varepsilon\to 0} \left(\limsup_{N\to \infty} \frac{\log \#\left(E, d_N,\varepsilon\right)}{N}\right). \]

Let $(X, T)$ and $(Y, S)$ be dynamical systems.
Let $\pi:X\to Y$ be an equivariant continuous map. Here ``equivariant’’ means $\pi\circ T = S\circ \pi$.
We often denote such a map by $\pi:(X, T)\to (Y, S)$ for clarifying the maps $T$ and $S$.
Let $d$ be a metric on $X$.
We define the \textbf{topological conditional entropy} of $\pi$ by 
\[  \htop(\pi, T)  =  \lim_{\varepsilon\to 0} \left(\lim_{N\to \infty} 
      \frac{\sup_{y\in Y}\log \#\left(\pi^{-1}(y), d_N, \varepsilon\right)}{N}\right). \]
This definition is due to \cite[Lemma 6.8.2]{Downarowicz}. 
It coincides with the supremum of \textit{fiber-wise entropy} \cite[Theorem 6.8.3]{Downarowicz}:
\[  \htop(\pi, T) = \sup_{y\in Y} \htop\left(\pi^{-1}(y), T\right). \]
In the book of Downarowicz \cite{Downarowicz} the topological conditional entropy $\htop(\pi, T)$
is denoted by $\mathbf{h}(T|S)$.

Bowen \cite[Theorem 17]{Bowen71} proved that for any equivariant continuous map 
$\pi:(X, T)\to (Y, S)$ we have
\begin{equation} \label{eq: Bowen inequality}
   \htop(X, T) \leq \htop(Y, S) + \htop(\pi, T).
\end{equation}
This inequality has a clear intuitive meaning; if $(X, T)$ is “larger than” $(Y, S)$ then some fiber of $\pi$ must be “large”.

We are interested in the case that both $\htop(X, T)$ and $\htop(Y, S)$ are infinite.
In this case Bowen’s inequality (\ref{eq: Bowen inequality}) provides no information about the fibers of $\pi$.
Typical examples of such infinite entropy systems are given by \textit{shifts on infinite dimensional cubes}.
Let $a$ be a natural number, and let $[0,1]^a$ be the $a$-dimensional cube.
We consider the bi-infinite product of $[0,1]^a$:
\[  \left([0,1]^a\right)^{\mathbb{Z}} = \cdots \times [0,1]^a\times [0,1]^a\times [0,1]^a\times \cdots. \]
We define the shift $\sigma: \left([0,1]^a\right)^{\mathbb{Z}}\to \left([0,1]^a\right)^{\mathbb{Z}}$ by 
\[  \sigma\left((x_n)_{n\in \mathbb{Z}}\right) = (x_{n+1})_{n\in \mathbb{Z}}, \quad (\text{each }x_n\in [0,1]^a). \]
The pair $\left(\left([0,1]^a\right)^{\mathbb{Z}}, \sigma\right)$ is a dynamical system which has infinite topological entropy.
This system is called the \textbf{full-shift on the alphabet $[0,1]^a$}.
(Later we will often consider two full-shifts on different alphabets simultaneously. 
But we use the same letter $\sigma$ for shift maps of both systems. Hopefully this will not cause any confusion.)

The following theorem is a starting point of our investigation.
\begin{theorem} \label{theorem: topological conditional entropy between two full shifts}
Let $a>b$ be two natural numbers. Let $\left(\left([0,1]^a\right)^{\mathbb{Z}}, \sigma\right)$ and 
$\left(\left([0,1]^b\right)^{\mathbb{Z}}, \sigma\right)$ be the full-shifts on the alphabets $[0,1]^a$ and $[0,1]^b$ respectively.
Then for any equivariant continuous map $\pi:\left([0,1]^a\right)^{\mathbb{Z}}\to \left([0,1]^b\right)^{\mathbb{Z}}$ we have 
\[   \htop(\pi, \sigma) = +\infty. \]
\end{theorem}

This theorem has the same spirit as that of Bowen’s inequality (\ref{eq: Bowen inequality}):
Since we assumed $a>b$,
the shift $\left([0,1]^a\right)^{\mathbb{Z}}$ is “larger than” the shift $\left([0,1]^b\right)^{\mathbb{Z}}$, and hence
some fiber of $\pi:\left([0,1]^a\right)^{\mathbb{Z}}\to \left([0,1]^b\right)^{\mathbb{Z}}$ must be “large”.
However this result does not follow from Bowen’s inequality (\ref{eq: Bowen inequality}) because
both the shifts $\left([0,1]^a\right)^{\mathbb{Z}}$ and $\left([0,1]^b\right)^{\mathbb{Z}}$ have infinite topological 
entropy. 

Theorem \ref{theorem: topological conditional entropy between two full shifts} might look an innocent statement.
However we can add it a small twist so that it becomes more interesting:

\begin{theorem} \label{theorem: topological conditional entropy and minimal system}
Let $a>b$ be two natural numbers.
There exists a minimal closed subset $X$ of $\left(\left([0,1]^a\right)^{\mathbb{Z}}, \sigma\right)$ such that 
for any equivariant continuous map $\pi:\left(X, \sigma\right) \to \left(\left([0,1]^b\right)^{\mathbb{Z}}, \sigma\right)$ we have 
\[   \htop(\pi, \sigma) = +\infty. \]
\end{theorem}

Here “minimal” means that $\sigma(X) = X$ and that every orbit of $(X, \sigma)$ is dense in $X$.

Letting $a=2$ and $b=1$ in the theorem, we get a minimal closed subset $X$ of the shift $\left([0,1]^2\right)^{\mathbb{Z}}$ such that 
for any equivariant continuous map $\pi: X\to [0,1]^{\mathbb{Z}}$ we have $\htop(\pi, \sigma) = +\infty$.
In particular, \textit{$(X, \sigma)$ does not embed in the shift $[0,1]^{\mathbb{Z}}$}.
This provides an answer to Auslander’s question in \cite[p. 193]{Auslander} from a new angle;
Auslander asked whether there exists a minimal dynamical system which does not embed in the shift $[0,1]^{\mathbb{Z}}$.
This problem was first solved in the celebrated paper of Lindenstrauss--Weiss \cite{Lindenstrauss--Weiss} by using mean dimension theory.
We briefly review their work below.

Mean dimension is a topological invariant of dynamical systems introduced by Gromov \cite{Gromov99}.
Let $(X, T)$ be a dynamical system with a metric $d$ on $X$.
For $\varepsilon>0$ and a natural number $N$ we define $\widim_\varepsilon(X, d_N)$ as the minimum number of $n$ such that 
there exist an $n$-dimensional finite simplicial complex $P$ and an $\varepsilon$-embedding 
$f:(X, d_N)\to P$. 
Here $f$ is called an \textbf{$\varepsilon$-embedding} if we have $\diam(f^{-1}(p), d_N) < \varepsilon$ for all $p\in P$.  
We define the \textbf{mean dimension} of $(X, T)$ by 
\[  \mdim(X, T) := \lim_{\varepsilon\to 0} \left(\lim_{N\to \infty} \frac{\widim_\varepsilon(X, d_N)}{N}\right). \]
This is a topological invariant, i.e. independent of the choice of $d$.
For a closed subset $E\subset X$ we also define 
\[  \mdim(E, T) := \lim_{\varepsilon\to 0} \left(\limsup_{N\to \infty} \frac{\widim_\varepsilon(E, d_N)}{N}\right). \]

If a dynamical system $(X, T)$ embeds in another $(Y, S)$ then we have $\mdim(X, T)\leq \mdim(Y, S)$.
The full-shift $\left(\left([0,1]^a\right)^{\mathbb{Z}}, \sigma\right)$ has mean dimension $a$.
Lindenstrauss-Weiss \cite[Proposition 3.5]{Lindenstrauss--Weiss}
constructed a minimal closed subset $X\subset \left([0,1]^2\right)^{\mathbb{Z}}$ such that 
$(X, \sigma)$ has mean dimension strictly larger than one.
Then such $X$ does not embed in $[0,1]^{\mathbb{Z}}$ because $[0,1]^{\mathbb{Z}}$ has mean dimension one.
So this provides a negative answer to Auslander’s question.

Notice that the statement “\textit{$X$ does not embed in $[0,1]^{\mathbb{Z}}$}” means that every
equivariant continuous map $\pi:X\to [0,1]^{\mathbb{Z}}$ has some fiber of cardinality strictly larger than one.
Theorem \ref{theorem: topological conditional entropy and minimal system} above provides a substantially stronger conclusion that 
some fiber of $\pi$ must have infinite entropy.
Indeed the construction of the minimal closed subset $X$ in Theorem \ref{theorem: topological conditional entropy and minimal system}
is the same as that of Lindenstrauss-Weiss \cite[Proposition 3.5]{Lindenstrauss--Weiss}.
But we can deduce a stronger conclusion with the aid of the waist inequality.

\subsection{Conditional metric mean dimension}  \label{subsection: conditional metric mean dimension}

As we already saw at the end of the last subsection, our study has a connection to mean dimension theory.
It turns out that the most fruitful framework is provided by a conditional version of \textit{metric mean dimension}.
Metric mean dimension is a notion introduced by Lindenstrauss--Weiss \cite[\S 4]{Lindenstrauss--Weiss}
for better understanding the relation between mean dimension and entropy.
Its conditional version was introduced by Liang \cite[Definition 4.1]{Liang}.

Let $(X, T)$ be a dynamical system with a metric $d$ on $X$.
We define the \textbf{upper and lower metric mean dimensions} of $(X, T, d)$ by 
\begin{align*}
   \umdimm(X, T, d) 
  & = \limsup_{\varepsilon\to 0} 
   \left(\lim_{N\to \infty} \frac{\log \#\left(X, d_N, \varepsilon\right)}{N \log(1/\varepsilon)}\right), \\
   \lmdimm(X, T, d)
  & =  \liminf_{\varepsilon\to 0} 
   \left(\lim_{N\to \infty} \frac{\log \#\left(X, d_N, \varepsilon\right)}{N \log(1/\varepsilon)}\right).
\end{align*}   
These are metric-dependent quantities.
When the upper and lower limits coincide, we denote the common value by 
$\mdimm(X, T, d)$.
If $\lmdimm(X, T, d) >0$ then $\htop(X, T)  = +\infty$.
So, in a sense, metric mean dimension evaluates \textit{how infinite topological entropy is}. 

Lindenstrauss--Weiss \cite[Theorem 4.2]{Lindenstrauss--Weiss} proved 
\[  \mdim(X, T) \leq \lmdimm(X, T, d). \]
See also the paper of Lindenstrauss \cite{Lindenstrauss} for a deeper study of metric mean dimension.

Let $(X, T)$ and $(Y, S)$ be dynamical systems, and let $\pi:X\to Y$ be an equivariant continuous map .
Let $d$ be a metric on $X$.
We define the \textbf{upper and lower conditional metric mean dimensions} of $\pi$ by 
\begin{align*}
   \umdimm(\pi, T, d) 
  & = \limsup_{\varepsilon\to 0} 
   \left(\lim_{N\to \infty} \frac{\sup_{y\in Y} \log \#\left(\pi^{-1}(y), d_N, \varepsilon\right)}{N \log(1/\varepsilon)}\right), \\
   \lmdimm(\pi, T, d)
  & =  \liminf_{\varepsilon\to 0} 
   \left(\lim_{N\to \infty} \frac{\sup_{y\in Y} \log \#\left(\pi^{-1}(y), d_N, \varepsilon\right)}{N \log(1/\varepsilon)}\right).
\end{align*} 
These are also metric-dependent quantities.
When the upper and lower limits coincide, we denote the common value by 
$\mdimm(\pi, T, d)$.
If $\lmdimm(\pi, T, d) >0$ then $\htop(\pi, T) = +\infty$.
So conditional metric mean dimension evaluates \textit{how infinite topological conditional entropy is}.

Let $a$ be a natural number. For $x = (x_1,\dots, x_a)\in \mathbb{R}^a$ we set 
\[   \norm{x}_\infty := \max_{1\leq i \leq a} |x_i|. \]

Let $\left([0,1]^a\right)^{\mathbb{Z}}$ be the full-shift on the alphabet $[0,1]^a$.
We define a metric $D$ on it by 
\begin{equation}  \label{eq: metric on full shift}
  D\left(\left(x_n\right)_{n\in \mathbb{Z}}, \left(y_n\right)_{n\in \mathbb{Z}}\right) 
    = \sum_{n\in \mathbb{Z}} 2^{-|n|} \norm{x_n-y_n}_\infty, \quad 
    (x_n, y_n\in [0,1]^a).
\end{equation}

\begin{theorem}  \label{theorem: conditional metric mean dimension between two full shifts}
Let $a>b$ be two natural numbers. Let $\left(\left([0,1]^a\right)^{\mathbb{Z}}, \sigma\right)$ and 
$\left(\left([0,1]^b\right)^{\mathbb{Z}}, \sigma\right)$ be the full-shifts on the alphabets $[0,1]^a$ and $[0,1]^b$ respectively.
Then for any equivariant continuous map $\pi:\left([0,1]^a\right)^{\mathbb{Z}}\to \left([0,1]^b\right)^{\mathbb{Z}}$ we have 
\[   \lmdimm(\pi, \sigma, D) \geq a - b. \]
Here $D$ is the metric defined by (\ref{eq: metric on full shift}).
\end{theorem}

Theorem \ref{theorem: topological conditional entropy between two full shifts} in \S \ref{subsection: topological conditional entropy}
immediately follows from this theorem because $\lmdimm(\pi, \sigma, D) > 0$ implies 
$\htop(\pi, \sigma) = +\infty$.

\begin{example}
The estimate in Theorem \ref{theorem: conditional metric mean dimension between two full shifts} is sharp:
For $a>b$, let $f:[0,1]^a\to [0,1]^b$ be the natural projection to the first $b$ coordinates.
We define $\pi:\left([0,1]^a\right)^{\mathbb{Z}}\to \left([0,1]^b\right)^{\mathbb{Z}}$ by 
\[  \pi\left((x_n)_{n\in \mathbb{Z}}\right) = \left(f(x_n)\right)_{n\in \mathbb{Z}}. \]
Then 
\[   \mdimm(\pi, \sigma, D) = a-b. \]
\end{example}

We can also prove a statement similar to Theorem \ref{theorem: conditional metric mean dimension between two full shifts} 
for maps from $\left([0,1]^a\right)^{\mathbb{Z}}$ to \textit{arbitrary} dynamical systems:

\begin{theorem}  \label{theorem: conditional metric mean dimension and mean dimension}
Let $a$ be a natural number and let $(Y, S)$ be a dynamical system.
For any equivariant continuous map $\pi:\left(\left([0,1]^a\right)^{\mathbb{Z}}, \sigma\right) \to (Y, S)$ we have 
\begin{equation} \label{eq: conditional metric mean dimension and mean dimension}
    \lmdimm(\pi, \sigma, D) \geq a - 2\mdim(Y, S). 
\end{equation}    
Here $\mdim(Y, S)$ is the mean dimension of $(Y, S)$.
\end{theorem}

The factor $2$ in the right-hand side of (\ref{eq: conditional metric mean dimension and mean dimension}) looks a bit unsatisfactory.
So we propose:

\begin{problem} \label{problem: reduce the factor two to one}
In the setting of Theorem \ref{theorem: conditional metric mean dimension and mean dimension}, can one conclude 
the following stronger inequality?
\[   \lmdimm(\pi, \sigma, D) \geq a - \mdim(Y, S).  \]
See also Remark \ref{remark: proof of the main technical theorem} (2).
\end{problem}

Here the target system $(Y, S)$ is arbitrary, but the domain is only $\left([0,1]^a\right)^{\mathbb{Z}}$.
It seems difficult for our current technology to try a similar question for maps $\pi:(X, T)\to (Y, S)$
with arbitrary dynamical systems $(X, T)$.
Nevertheless we would like to propose the following problem because the statement looks beautiful (if it is true).

\begin{problem}
Let $(X, T)$ and $(Y, S)$ be dynamical systems with 
\[   \mdim(X, T) > \mdim(Y, S).  \]
Let $\pi: (X, T) \to (Y, S)$ be an equivariant continuous map.
Can one always conclude that the topological conditional entropy of $\pi$ is infinite?
When $\mdim(Y, S) = 0$, this can be proved by using \cite[Theorem 1.5]{Tsukamoto analogue of Hurewicz}.
But when $\mdim(Y, S)>0$, the problem is widely open.
\end{problem}

Although it is currently difficult to study maps from general dynamical systems,
our method works well for some subsystems of $\left([0,1]^a\right)^{\mathbb{Z}}$.
We can prove the following theorem.

\begin{theorem} \label{theorem: conditional metric mean dimension and minimal system}
Let $a$ be a natural number.
Let $s$ be a real number with $0\leq s < a$. 
There exists a minimal closed subset $X$ of $\left(\left([0,1]^a\right)^{\mathbb{Z}}, \sigma\right)$ satisfying the following two conditions.
 \begin{enumerate}
    \item $\mdim(X, \sigma) = \mdimm\left(X, \sigma, D\right) = s$.
    \item For any natural number $b$ and any equivariant continuous map 
    $\pi:(X, \sigma) \to \left(\left([0,1]^b\right)^{\mathbb{Z}}, \sigma\right)$ we have 
    \[   \lmdimm(\pi, \sigma, D) \geq s-b. \]
 \end{enumerate}
Here $D$ is a metric defined in (\ref{eq: metric on full shift}).
\end{theorem}

When $s>b$, the inequality $\lmdimm(\pi, \sigma, D) \geq s-b >0$ implies $\htop(\pi, \sigma) = +\infty$.
Therefore Theorem \ref{theorem: topological conditional entropy and minimal system} in \S \ref{subsection: topological conditional entropy}
is an immediate corollary of Theorem \ref{theorem: conditional metric mean dimension and minimal system}.

\subsection{On an analogue of the Hurewicz theorem}  \label{subsection: on an analogue of the Hurewicz theorem}

The Hurewicz theorem \cite[p.91, Theorem VI 7]{Hurewicz--Wallman} in classical dimension theory asserts that
for any continuous map $f:X\to Y$ between compact metrizable spaces we have 
\begin{equation} \label{eq: Hurewicz}
   \dim X \leq \dim Y + \sup_{y\in Y} \dim f^{-1}(y). 
\end{equation}   
Here $\dim(\cdot)$ denotes topological dimension\footnote{Let $d$ be a metric on $X$.
We define the \textbf{topological dimension} by $\dim X = \lim_{\varepsilon\to 0} \widim_\varepsilon(X, d)$.
The value of $\dim X$ is independent of the choice of $d$}.

Let $d$ be a metric on $X$. 
We define the \textbf{upper and lower Minkowski dimensions} of $(X, d)$ by 
\[ \udimm(X, d) = \limsup_{\varepsilon\to 0} \frac{\log \#\left(X, d, \varepsilon\right)}{\log(1/\varepsilon)}, \quad
 \underline{\dim}_{\mathrm{M}}(X, d) 
   = \liminf_{\varepsilon\to 0} \frac{\log \#\left(X, d, \varepsilon\right)}{\log(1/\varepsilon)}. \]
When the upper and lower limits coincide, we denote the common value by 
$\dim_{\mathrm{M}}(X, d)$.
The topological dimension is always bounded by the lower Minkowski dimension.
\[  \dim X \leq \underline{\dim}_{\mathrm{M}}(X, d). \]
So $\dim f^{-1}(y) \leq \ldimm\left(f^{-1}(y), d\right)$, and hence we also have 
\begin{equation} \label{eq: Hurewicz and Minkowski}
   \dim X \leq \dim Y + \sup_{y\in Y} \ldimm\left(f^{-1}(y), d\right)
\end{equation}
for every continuous map $f:X\to Y$.

Theorems \ref{theorem: conditional metric mean dimension between two full shifts} and 
\ref{theorem: conditional metric mean dimension and mean dimension} in \S 
\ref{subsection: conditional metric mean dimension} are dynamical analogues of 
the inequality (\ref{eq: Hurewicz and Minkowski}).
Conditional metric mean dimension is an analogue of the term $\sup_{y\in Y} \ldimm\left(f^{-1}(y),d\right)$.
Probably it is also natural to ask whether or not an analogue of the original Hurewicz theorem (\ref{eq: Hurewicz}) holds for mean dimension.
Indeed this question was already proposed by the second-named author \cite{Tsukamoto moduli space of Brody curves} 
more than ten years ago and studied in detail by a recent paper \cite{Tsukamoto analogue of Hurewicz}.
In this subsection we explain the relation between \cite{Tsukamoto analogue of Hurewicz} and the current paper.

The paper \cite[Theorem 3.9, Remark 3.10]{Tsukamoto analogue of Hurewicz} proved the next theorem.

\begin{theorem} \label{theorem: analogue of Hurewicz does not hold}
Let $a$ be a natural number. There exists a zero-dimensional\footnote{“Zero-dimensional” means that the 
topological dimension of $Z$ is zero.} dynamical system $(Z, R)$ for which the following 
statement holds true:
For any $\delta>0$ there exist a dynamical system $(Y, S)$ and an equivariant continuous map 
\[  \pi:\left(\left([0,1]^a\right)^{\mathbb{Z}}\times Z, \sigma\times R\right) \to (Y, S)  \]
satisfying 
\[  \mdim(Y, S) < \delta, \quad  \sup_{y\in Y} \mdim\left(\pi^{-1}(y), \sigma\times R\right) = 0. \]
Here $\mdim\left(\pi^{-1}(y), \sigma\times R\right)$ is the mean dimension of the fiber 
$\pi^{-1}(y) \subset \left([0,1]^a\right)^{\mathbb{Z}}\times Z$.
\end{theorem}

In this statement, the system $\left(\left([0,1]^a\right)^{\mathbb{Z}}\times Z, \sigma\times R\right)$ has mean dimension $a$
while $\mdim(Y, S) < \delta$ may be arbitrarily small.
Nonetheless the map $\pi$ has no mean dimension in the fiber direction.
This implies that a direct analogue of the Hurewicz theorem (\ref{eq: Hurewicz}) does not hold for mean dimension.
In other words, mean dimension of the fibers $\mdim\left(\pi^{-1}(y), \sigma\times R\right)$
does not properly evaluate the complexity of the map $\pi$.

On the other hand, we can prove the following statement for conditional metric mean dimension:

\begin{proposition}  \label{prop: analogue of metric Hurewicz}
Let $a$ be a natural number. Let $(Z, R)$ be a dynamical system with a metric $\rho$ on $Z$.
We define a metric $d$ on the product $\left([0,1]^a\right)^{\mathbb{Z}}\times Z$ by 
\[  d\left((x, z), (x^\prime, z^\prime)\right) = D(x, x^\prime) + \rho(z,z^\prime), \quad  \]
where $D$ is the metric on $\left([0,1]^a\right)^{\mathbb{Z}}$ defined by (\ref{eq: metric on full shift}).
Then for any dynamical system $(Y, S)$ and any equivariant continuous map 
$\pi: \left(\left([0,1]^a\right)^{\mathbb{Z}}\times Z, \sigma\times R\right) \to (Y, S)$ we have 
\[  \lmdimm\left(\pi, \sigma\times R, d\right) \geq a - 2\mdim(Y, S). \]
\end{proposition}

Namely, we can detect a reasonable amount of complexity of the map $\pi$ by using conditional metric mean dimension.
This gives the authors more confidence that the conditional metric mean dimension is a useful tool for studying
maps between dynamical systems.

\begin{remark}
In the setting of Proposition \ref{prop: analogue of metric Hurewicz}, we can also prove the following statement:
Define an equivariant continuous map 
$\Pi:  \left(\left([0,1]^a\right)^{\mathbb{Z}}\times Z, \sigma\times R\right) \to (Y\times Z, S\times R)$ by 
\[  \Pi(x, z) = \left(\pi(x,z), z\right). \]
Then we have 
\[  \lmdimm\left(\Pi, \sigma \times R, d\right) \geq a - 2\mdim(Y, S). \]
This statement is a bit stronger (and arguably more natural) than Proposition \ref{prop: analogue of metric Hurewicz}.
But the proof is more or less the same.
So we omit the details of the proof.
\end{remark}

\subsection{Organization of the paper}

We review the waist inequality in the next section.
We prepare a simple fact on conditional metric mean dimension in \S \ref{section: a preparation on conditional metric mean dimension}.
We prove the main technical theorem in \S \ref{section: main technical theorem}.
As its corollaries, we prove Theorems \ref{theorem: conditional metric mean dimension between two full shifts}, 
\ref{theorem: conditional metric mean dimension and mean dimension} and Proposition \ref{prop: analogue of metric Hurewicz} 
also in \S \ref{section: main technical theorem}.
We prove Theorem \ref{theorem: conditional metric mean dimension and minimal system} 
in \S \ref{section: construction of a minimal subset}.

\section{Waist inequality}   \label{section: waist inequality}

In this section we review the waist inequality.
Readers can find much more information in an excellent survey of Guth \cite{Guth}.
Let $S^n = \{x\in \mathbb{R}^{n+1}\mid |x| = 1\}$ be the $n$-dimensional sphere.
Here $|x|$ is the Euclidean norm.
We naturally consider $S^0 \subset S^1 \subset S^2\subset \dots \subset S^n$.
For $A\subset S^n$ and $r>0$ we denote by $A+r$ the $r$-neighborhood of $A$ in the spherical metric.

The following inequality is the original waist inequality proved by Gromov \cite{Gromov03}.

\begin{theorem}[Gromov, 2003]
Let $n\geq m$ be two natural numbers.
For any continuous map $f:S^n\to \mathbb{R}^m$ there exists $t\in \mathbb{R}^m$ such that for every $r>0$ we have 
\begin{equation} \label{eq: Gromov waist inequality}
   \vol\left(f^{-1}(t)+r\right) \geq \vol\left(S^{n-m} + r\right). 
\end{equation}   
Here $\vol(\cdot)$ denotes the $n$-dimensional volume on the sphere $S^n$.
\end{theorem}

For maps $f:S^n\to \mathbb{R}$, this statement follows from the isoperimetric inequality on the sphere.
(See \cite[Appendix]{Alpert--Guth} for the details.)
For equidimensional maps $f:S^n\to \mathbb{R}^n$, if we additionally assume that every fiber of $f$ has cardinality at most two, 
then the above inequality (\ref{eq: Gromov waist inequality}) for $r = \frac{\pi}{2}$ implies that 
$f^{-1}(t)$ consists of a pair of antipodal points. 
This is the statement of the Borsuk--Ulam theorem for such maps $f$.
Therefore the waist inequality has close connections to both the isoperimetric inequality and the Borsuk--Ulam theorem.
Its proof is a wonderful mixture of geometry (measure theory and convex geometry) and algebraic topology.
Gromov’s original paper \cite{Gromov03} is probably not very easily readable.
Memarian’s paper \cite{Memarian} is helpful for understanding the details of the proof of the waist inequality.

As we already mentioned in the introduction, the ideas of the waist inequality have connections to several subjects.
Here we review a connection to convex geometry developed by Klartag \cite{Klartag17}.
This will be crucial for us.

For subsets $A$ and $B$ of $\mathbb{R}^n$ we define 
$A+B = \{x+y\mid x\in A, y\in B\}$.
For $r>0$ we set $r A = \{rx\mid x\in A\}$.
A \textbf{convex body} of $\mathbb{R}^n$ is a compact convex subset of $\mathbb{R}^n$ with a non-empty interior.
A convex body $K$ is said to be \textbf{centrally-symmetric} if $-K=K$, where 
$-K := \{-x\mid x\in K\}$.
A function $\varphi:\mathbb{R}^n\to [0, +\infty)$ is said to be \textbf{log-concave} if for any $x, y\in \mathbb{R}^n$ and $0<\lambda<1$
\[   \varphi\left(\lambda x + (1-\lambda)y\right) \geq \varphi(x)^\lambda \varphi(y)^{1-\lambda}. \]

Klartag \cite{Klartag17} proved several interesting theorems connecting the waist inequality to convex geometry.
In particular he proved, among many other things, the following version of the waist inequality \cite[Theorem 5.7]{Klartag17}.

\begin{theorem}[Klartag, 2017]
Let $K\subset \mathbb{R}^n$ be a centrally-symmetric convex body.
Let $1\leq m \leq n$. Let $\mu$ be a probability measure supported in $K$ which has a log-concave density function 
with respect to the Lebesgue measure.
Then for any continuous map $f:K\to \mathbb{R}^m$ there exists $t\in \mathbb{R}^m$ satisfying 
\[  \mu\left(f^{-1}(t)+rK\right) \geq \left(\frac{r}{2+r}\right)^m \> \text{ for all $0<r < 1$}. \]
\end{theorem}

For $A\subset \mathbb{R}^n$ and $r>0$ we define $A+_\infty r$ as the 
$r$-neighborhood of $A$ with respect to the $\ell^\infty$-norm:
\[  A+_\infty r = \{x\in \mathbb{R}^n\mid \exists a\in A: \norm{x-a}_\infty \leq r\}. \]

\begin{corollary}  \label{corollary:  Klartag waist inequality}
Let $1\leq m \leq n$. Let $\mu$ be the Lebesgue measure restricted to the unit cube $[0,1]^n$.
Then for any continuous map $f:[0,1]^n\to \mathbb{R}^m$ there exists $t\in \mathbb{R}^m$ satisfying 
\[   \mu\left(f^{-1}(t) +_\infty r \right) \geq \left(\frac{r}{2}\right)^m \> \text{ for all $0 < r < \frac{1}{2}$}. \]
\end{corollary}

\begin{proof}
Let $K = \left[-\frac{1}{2}, \frac{1}{2}\right]^n$.
This is a centrally-symmetric convex body.
Let $\nu$ be the Lebesgue measure restricted to $K$.
The density of $\nu$ is the characteristic function of $K$, which is obviously log-concave.
We define $g:K\to \mathbb{R}^m$ by 
\[   g(x)  = f\left(x-\left(\frac{1}{2}, \frac{1}{2}, \dots, \frac{1}{2}\right)\right).   \]
By Klartag’s waist inequality, there exists $t\in \mathbb{R}^n$ satisfying 
\[  \nu\left(g^{-1}(t) + r K\right) \geq \left(\frac{r}{2+r}\right)^m \> \text{ for all $0<r < 1$}. \]
We have 
\[   g^{-1}(t) +_\infty r = g^{-1}(t) + 2r K. \]
Hence, for $0<r<\frac{1}{2}$
\[
   \nu\left(g^{-1}(t)+_\infty r\right)  = \nu\left(g^{-1}(t) + 2r K\right)  
    \geq \left(\frac{2r}{2+2r}\right)^m  \geq \left(\frac{r}{2}\right)^m.
\]
We have 
\[    f^{-1}(t) +_\infty r = \left(g^{-1}(t)+_\infty r\right) + \left(\frac{1}{2}, \frac{1}{2}, \dots, \frac{1}{2}\right). \]
Thus, for $0 < r < \frac{1}{2}$
\[  \mu\left(f^{-1}(t) +_\infty r\right) = \nu\left(g^{-1}(t)+_\infty r\right)  \geq  \left(\frac{r}{2}\right)^m. \]
\end{proof}

For $r>0$ and a compact subset $A\subset \mathbb{R}^n$ we define $\#\left(A, \norm{\cdot}_\infty, r\right)$
as the minimum number of $N$ such that there exists an open cover $\{U_1, U_2, \dots, U_N\}$ of $A$ satisfying 
$\diam\left(U_i, \norm{\cdot}_\infty\right) < r$ for all $1\leq i \leq r$.
Here $\diam\left(U_i, \norm{\cdot}_\infty\right) = \sup_{x, y\in U_i} \norm{x-y}_\infty$.
When $A$ is empty, we define $\#\left(A, \norm{\cdot}_\infty, r\right) = 0$.

\begin{corollary} \label{corollary: waist inequality for covering number}
Let $1\leq m \leq n$. 
For any continuous map $f:[0,1]^n\to \mathbb{R}^m$ there exists $t\in \mathbb{R}^m$ satisfying 
\[   \#\left(f^{-1}(t), \norm{\cdot}_\infty, r\right) \geq \frac{1}{8^n} \left(\frac{1}{r}\right)^{n-m} \> \text{ for all $0 < r < \frac{1}{2}$}. \]
\end{corollary}

\begin{proof}
Let $\mu$ be the Lebesgue measure restricted to $[0,1]^n$.
By Corollary \ref{corollary:  Klartag waist inequality} there exists $t\in \mathbb{R}^m$ satisfying 
\[  \mu\left(f^{-1}(t) +_\infty r\right) \geq \left(\frac{r}{2}\right)^m \> \text{ for all $0 < r < \frac{1}{2}$}. \]
Let $0 < r < \frac{1}{2}$.
Set $N = \#\left(f^{-1}(t), \norm{\cdot}_\infty, r\right)$.
There exists an open cover $f^{-1}(t) \subset U_1\cup U_2\cup \dots \cup U_N$ satisfying 
$\diam\left(U_i, \norm{\cdot}_\infty\right) < r$ for all $1\leq i \leq N$.
Pick a point $x_i\in U_i$ for each $U_i$.
We have 
\[  U_i \subset x_i + [-r, r]^n := \left\{y\in \mathbb{R}^n\mid \norm{y-x_i}_\infty \leq r\right\}.  \]
Hence 
\[   f^{-1}(t) \subset \bigcup_{i=1}^N \left(x_i + [-r, r]^n\right). \]
We have 
\[   f^{-1}(t) +_\infty r = \bigcup_{i=1}^N \left(x_i + [-2r, 2r]^n\right). \]
Then 
\[
  \left(\frac{r}{2}\right)^m \leq \mu\left(f^{-1}(t) +_\infty r\right) 
   \leq \sum_{i=1}^N \mu\left(x_i + [-2r, 2r]^n\right)  \leq N (4r)^n.
\]
Therefore 
\[   N \geq (4r)^{-n} \cdot  \left(\frac{r}{2}\right)^m \geq  \frac{1}{8^n}\cdot \left(\frac{1}{r}\right)^{n-m}. \] 
\end{proof}

For a convenience of readers, here we recall basic notions used in the definition of mean dimension.
Let $\varepsilon$ be a positive number, and 
let $(X, d)$ be a compact metric space.
A continuous map $f:X\to Y$ from $X$ to a topological space $Y$ is called an $\varepsilon$-embedding
if we have $\diam f^{-1}(y) < \varepsilon$ for all $y\in Y$.
We define $\widim_\varepsilon(X, d)$ as the minimum number $n$ such that there exist an $n$-dimensional finite simplicial complex $P$ 
and an $\varepsilon$-embedding $f:X\to P$.

We will need the following lemma later.
This was proved in \cite[Lemma 3.2]{Lindenstrauss--Weiss}.
\begin{lemma}  \label{lemma: Lebesgue}
For $0<\varepsilon<1$ we have 
\[   \widim_\varepsilon\left([0,1]^n, \norm{\cdot}_\infty\right) = n. \]
Here $\norm{\cdot}_\infty$ is the metric defined by the $\ell^\infty$-norm.
\end{lemma}

This lemma is a basis of the fact that the full-shift $\left([0,1]^a\right)^{\mathbb{Z}}$ has mean dimension $a$.

\section{A preliminary on conditional metric mean dimension}  \label{section: a preparation on conditional metric mean dimension}

In this section we prepare a simple fact about conditional metric mean dimension.

Let $(X, T)$ be a dynamical system with a metric $d$.
Recall that, for each $N\geq 1$, we defined a new metric $d_N$ on $X$ by 
$d_N(x_1, x_2) = \max_{0\leq n < N} d\left(T^n x_1, T^n x_2\right)$.
Let $(Y, S)$ be another dynamical system with a metric $d^\prime$.
Let $N$ be a natural number and let $r$ be a positive real number.
For $y\in Y$ we set 
\[  B_r\left(y, d^\prime_N\right)  = \left\{x\in Y\mid d^\prime_N(x, y) \leq r\right\}. \]
Let $\pi:(X, T)\to (Y, S)$ be an equivariant continuous map.
Recall that we defined the upper and lower conditional metric mean dimensions of $\pi$ by 

\begin{align*}
   \umdimm(\pi, T, d) 
  & = \limsup_{\varepsilon\to 0} 
   \left(\lim_{N\to \infty} \frac{\sup_{y\in Y} \log \#\left(\pi^{-1}(y), d_N, \varepsilon\right)}{N \log(1/\varepsilon)}\right), \\
   \lmdimm(\pi, T, d)
  & =  \liminf_{\varepsilon\to 0} 
   \left(\lim_{N\to \infty} \frac{\sup_{y\in Y} \log \#\left(\pi^{-1}(y), d_N, \varepsilon\right)}{N \log(1/\varepsilon)}\right).
\end{align*} 
In this definition, the limit with respect to $N$ exists because 
the quantity 
\[  \sup_{y\in Y} \log \#\left(\pi^{-1}(y), d_N, \varepsilon\right) \]
is subadditive in $N$, as we will see below:

\begin{lemma}  \label{lemma: subadditivity}
In the above setting, we have the following two statements.
  \begin{enumerate}
      \item The quantity $\sup_{y\in Y} \log \#\left(\pi^{-1}(y), d_N, \varepsilon\right)$ is subadditive in $N$.
       \item Let $\varepsilon$ and $\delta$ be positive numbers. Then the quantity 
       \[   a_N := \sup_{y\in Y} \log \#\left(\pi^{-1}\left(B_\delta(y, d^\prime_N)\right), d_N, \varepsilon\right) \]
       is subadditive in $N$.       
  \end{enumerate}
\end{lemma}

\begin{proof}
  We only prove (2). The proof of (1) is similar.
Let $N_1$ and $N_2$ be natural numbers. We have 
\[  B_\delta\left(y, d^\prime_{N_1+N_2}\right) 
     = B_\delta\left(y, d^\prime_{N_1}\right) \cap S^{-N_1} B_\delta\left(S^{N_1}y, d^\prime_{N_2}\right). \]
Then 
\begin{align*}
    \pi^{-1}\left(B_\delta\left(y, d^\prime_{N_1+N_2}\right)\right) &
     = \pi^{-1}\left(B_\delta\left(y, d^\prime_{N_1}\right)\right) \cap \pi^{-1}\left(S^{-N_1} B_\delta\left(S^{N_1}y, d^\prime_{N_2}\right)\right) \\
    & = \pi^{-1}\left(B_\delta\left(y, d^\prime_{N_1}\right)\right) \cap T^{-N_1} \pi^{-1}\left(B_\delta\left(S^{N_1}y, d^\prime_{N_2}\right)\right).
\end{align*}
Hence 
\begin{align*}
    & \#\left( \pi^{-1}\left(B_\delta\left(y, d^\prime_{N_1+N_2}\right)\right), d_{N_1+N_2}, \varepsilon\right) \\
    &  \leq \#\left( \pi^{-1}\left(B_\delta\left(y, d^\prime_{N_1}\right)\right), d_{N_1}, \varepsilon\right) \cdot 
            \#\left(\pi^{-1}\left(B_\delta\left(S^{N_1}y, d^\prime_{N_2}\right)\right), d_{N_2}, \varepsilon\right) \\
    & \leq e^{a_{N_1}} \cdot e^{a_{N_2}}.        
\end{align*}            
Therefore we have $a_{N_1+N_2} \leq a_{N_1} + a_{N_2}$.
\end{proof}

\begin{lemma} \label{lemma: formula of conditional metric mean dimension}
 The upper and lower conditional metric mean dimensions of $\pi: (X, T) \to (Y, S)$ are given by
  \begin{align*}
      \umdimm(\pi, T, d) 
  & = \limsup_{\varepsilon\to 0} 
   \left\{\lim_{\delta\to 0} \left(\lim_{N\to \infty} 
   \frac{\sup_{y\in Y} \log \#\left(\pi^{-1}\left(B_\delta\left(y, d^\prime_N\right)\right), d_N, \varepsilon\right)}{N \log(1/\varepsilon)}\right)\right\}, \\
   \lmdimm(\pi, T, d)
  & =  \liminf_{\varepsilon\to 0} 
   \left\{\lim_{\delta\to 0} \left(\lim_{N\to \infty} 
   \frac{\sup_{y\in Y} \log \#\left(\pi^{-1}\left(B_\delta\left(y, d^\prime_N\right)\right), d_N, \varepsilon\right)}{N \log(1/\varepsilon)}\right)\right\}.
  \end{align*}
\end{lemma}

\begin{proof}
We only prove the formula for the lower conditional metric mean dimension.
The upper case is similar.
We set 
\[  \lmdimm(\pi, T, d)^\prime = 
     \liminf_{\varepsilon\to 0} 
   \left\{\lim_{\delta\to 0} \left(\lim_{N\to \infty} 
   \frac{\sup_{y\in Y} \log \#\left(\pi^{-1}\left(B_\delta\left(y, d^\prime_N\right)\right), d_N, \varepsilon\right)}{N \log(1/\varepsilon)}\right)\right\}. \]
It is obvious that $\lmdimm(\pi, T, d) \leq \lmdimm(\pi, T, d)^\prime$.
We will prove the reverse inequality.

Suppose $\lmdimm(\pi, T, d) < a$ for some positive number $a$.
We show $\lmdimm(\pi, T, d)^\prime \leq a$.
There exists a decreasing sequence $\varepsilon_1>\varepsilon_2>\varepsilon_3>\dots \to 0$ such that 
\[  \lim_{N\to \infty} \frac{\sup_{y\in Y} \log \#\left(\pi^{-1}(y), d_N, \varepsilon_k\right)}{N \log (1/\varepsilon_k)} < a.\]
There exists $N_k$ for each $k$ satisfying 
\[  \sup_{y\in Y} \#\left(\pi^{-1}(y), d_{N_k}, \varepsilon_k\right)  < \left(\frac{1}{\varepsilon_k}\right)^{a N_k}. \]
From the compactness of $Y$, we can find $\delta_k>0$ such that 
\[   \sup_{y\in Y} \#\left(\pi^{-1}\left(B_{\delta_k}(y, d^\prime_{N_k})\right), d_{N_k}, \varepsilon_k\right) 
     < \left(\frac{1}{\varepsilon_k}\right)^{a N_k}. \]
From the subadditivity (2) in Lemma \ref{lemma: subadditivity}
\[  \lim_{N\to \infty} \frac{\sup_{y\in Y} \log \#\left(\pi^{-1}\left(B_{\delta_k}(y, d^\prime_N)\right), d_N, \varepsilon_k\right)}{N}
      < a \log (1/\varepsilon_k). \]
Hence
\[  \lim_{\delta\to 0} \lim_{N\to \infty} \frac{\sup_{y\in Y} 
     \log \#\left(\pi^{-1}\left(B_{\delta}(y, d^\prime_N)\right), d_N, \varepsilon_k\right)}{N}
      < a \log (1/\varepsilon_k). \]
Thus we have $\lmdimm(\pi, T, d)^\prime \leq a$.
This has shown  $\lmdimm(\pi, T, d)^\prime \leq \lmdimm(\pi, T, d)$.
\end{proof}

\section{Main technical theorem}  \label{section: main technical theorem}

Here we formulate the main technical result of this paper.
All the theorems in \S \ref{section: introduction} are more or less its corollaries.
We prove Theorems \ref{theorem: conditional metric mean dimension between two full shifts}, 
\ref{theorem: conditional metric mean dimension and mean dimension} and Proposition \ref{prop: analogue of metric Hurewicz}
in this section.
Theorem \ref{theorem: conditional metric mean dimension and minimal system} will be proved in 
\S \ref{section: construction of a minimal subset}.

\begin{theorem}[Main Technical theorem]  \label{theorem: main technical theorem}
Let $s$ be a nonnegative real number.
Let $(X, T)$ be a dynamical system with a metric $d$ on $X$.
Suppose that there exist sequences of natural numbers $\{N_n\}_{n=1}^\infty$, $\{M_n\}_{n=1}^\infty$ and
continuous maps $\psi_n: [0,1]^{M_n}\to X$ $(n\geq 1)$ satisfying the following three conditions.
\begin{itemize}
  \item $N_n\to +\infty$ as $n\to \infty$.
  \item The ratio $M_n/N_n$ converges to $s$ as $n\to \infty$.
  \item For any $x, y\in [0,1]^{M_n}$ we have 
  \[  \norm{x-y}_\infty \leq  d_{N_n}\left(\psi_n(x), \psi_n(y)\right). \]
\end{itemize}
Then we have the following three conclusions.
  \begin{enumerate}
    \item $\mdim(X, T) \geq s$.
    \item For any natural number $b$ and any equivariant continuous map 
             $\pi:(X, T)\to \left(\left([0,1]^b\right)^{\mathbb{Z}}, \sigma\right)$ (where $\sigma$ is the shift map) we have 
            \[   \lmdimm\left(\pi, T, d\right) \geq s - b. \] 
    \item For any dynamical system $(Y, S)$ and any equivariant continuous map 
             $\pi:(X, T)\to (Y, S)$ we have 
             \[  \lmdimm\left(\pi, T, d\right) \geq s - 2\mdim(Y, S). \]      
  \end{enumerate}
\end{theorem}

\begin{proof}
(1) For $0<\varepsilon<1$ we have 
\[
   \widim_\varepsilon\left(X, d_{N_n}\right)   \geq \widim_\varepsilon\left([0,1]^{M_n}, \norm{\cdot}_\infty\right)  = M_n
   \quad \text{by Lemma \ref{lemma: Lebesgue}}.
\]
Then 
\[  \mdim\left(X, T\right) = \lim_{\varepsilon\to 0} \left(\lim_{n\to \infty} \frac{\widim_\varepsilon\left(X, d_{N_n}\right)}{N_n}\right) 
     \geq \lim_{n\to \infty} \frac{M_n}{N_n} = s. \]

(2) For two integers $\ell<m$ and a point $y = (y_n)_{n\in \mathbb{Z}}$ in $\left([0,1]^b\right)^{\mathbb{Z}}$ we denote 
\[    y|_\ell^m = \left(y_{\ell}, y_{\ell+1}, \dots, y_{m-1}, y_m\right) \in \left([0,1]^b\right)^{m-\ell+1}. \]
We take a metric $d^\prime$ on $\left([0,1]^b\right)^{\mathbb{Z}}$.
 Let $0<\varepsilon<1/2$ and $\delta>0$ be arbitrary. We take a natural number $m$ such that for 
 $y, z \in \left([0,1]^b\right)^{\mathbb{Z}}$
\[   y|_{-m}^m = z|_{-m}^m   \Longrightarrow  d^\prime(y, z) < \delta. \]
This implies that for any $N>0$ and $y, z \in \left([0,1]^b\right)^{\mathbb{Z}}$
\[   y|_{-m}^{N+m-1} = z|_{-m}^{N+m-1}   \Longrightarrow  d_N^\prime(y, z) < \delta. \]

We consider a map 
\[  \pi_n: [0,1]^{M_n} \to \left([0,1]^b\right)^{N_n+2m}, \quad 
     x \mapsto \pi\left(\psi_n(x)\right)|_{-m}^{N_n+m-1}.  \]
We apply the waist inequality to this map.
By Corollary \ref{corollary: waist inequality for covering number} we can find a point $t\in \left([0,1]^b\right)^{N_n+2m}$ satisfying 
\begin{equation}   \label{eq: consequence of waist inequality in main technical theorem} 
      \#\left(\pi_n^{-1}(t), \norm{\cdot}_\infty, \varepsilon\right) 
       \geq \frac{1}{8^{M_n}}\cdot \left(\frac{1}{\varepsilon}\right)^{M_n-b(N_n+2m)}. 
\end{equation}       
(When $M_n \geq b(N_n+2m)$, this is a direct consequence of Corollary \ref{corollary: waist inequality for covering number}.
When $M_n< b(N_n+2m)$, the right-hand side is less than one. So it is obviously true.)
Here, of course, the choice of $t$ depends on $n$. But we suppress its dependence on $n$ in our notation for simplicity.

Take $t^\prime\in \left([0,1]^b\right)^{\mathbb{Z}}$ with $t^\prime|_{-m}^{N_n+m-1} = t$.
We claim that 
\begin{equation} \label{eq: image of the fiber of pi_n}
    \psi_n\left(\pi_n^{-1}(t)\right) \subset \pi^{-1}\left(B_\delta\left(t^\prime, d^\prime_{N_n}\right)\right). 
\end{equation}    
Indeed, let $x\in \pi_n^{-1}(t)$. Then 
\[   \pi\left(\psi_n(x)\right)|_{-m}^{N_n+m-1} = \pi_n(x) = t = t^\prime|_{-m}^{N_n+m-1}. \]
By the choice of $m$ we have $d^\prime_{N_n}\left(\pi(\psi_n(x)), t^\prime\right) < \delta$.
Hence $\pi(\psi_n(x)) \in B_\delta\left(t^\prime, d^\prime_{N_n}\right)$ and 
$\psi_n(x) \in \pi^{-1}\left(B_\delta\left(t^\prime, d^\prime_{N_n}\right)\right)$.
This has proved (\ref{eq: image of the fiber of pi_n}).

By $\norm{x-y}_\infty \leq d_{N_n}\left(\psi_n(x), \psi_n(y)\right)$ $(x, y\in [0,1]^{M_n})$
\begin{align*}
   \#\left(\pi_n^{-1}(t), \norm{\cdot}_\infty, \varepsilon\right) &\leq \#\left(\psi_n\left(\pi_n^{-1}(t)\right), d_{N_n}, \varepsilon\right) \\
   & \leq \#\left(\pi^{-1}\left(B_\delta\left(t^\prime, d^\prime_{N_n}\right)\right), d_{N_n}, \varepsilon\right) \quad 
   \text{by (\ref{eq: image of the fiber of pi_n})}.
\end{align*}
By (\ref{eq: consequence of waist inequality in main technical theorem})
\[  \sup_{y\in \left([0,1]^b\right)^{\mathbb{Z}}} \#\left(\pi^{-1}\left(B_\delta\left(y, d^\prime_{N_n}\right)\right), d_{N_n}, \varepsilon\right) 
      \geq  \frac{1}{8^{M_n}}\cdot \left(\frac{1}{\varepsilon}\right)^{M_n-b(N_n+2m)}.  \]
Taking the logarithm, we have 
\[  \sup_{y\in \left([0,1]^b\right)^{\mathbb{Z}}} \log \#\left(\pi^{-1}\left(B_\delta\left(y, d^\prime_{N_n}\right)\right), d_{N_n}, \varepsilon\right) 
     \geq  -M_n \log 8 + (M_n-b N_n - 2bm) \log\left(\frac{1}{\varepsilon}\right). \]
We divide this by $N_n\log(1/\varepsilon)$ and let $n\to \infty$. By $\lim_{n\to \infty} M_n/N_n = s$
\[  \lim_{N\to \infty} 
\frac{\sup_{y\in \left([0,1]^b\right)^{\mathbb{Z}}} 
\log \#\left(\pi^{-1}\left(B_\delta\left(y, d^\prime_{N}\right)\right), d_{N}, \varepsilon\right)}{N\log(1/\varepsilon)}  
 \geq  s-b - \frac{s\log 8}{\log(1/\varepsilon)}. \]
This holds for arbitrary $0<\varepsilon<1/2$ and $\delta>0$.
So we can let $\delta\to 0$ and then $\varepsilon\to 0$.
By Lemma \ref{lemma: formula of conditional metric mean dimension}
\[    \lmdimm\left(\pi, T, d\right) \geq s-b. \]

(3)
Let $d^\prime$ be a metric on $Y$.
Let $0<\varepsilon<1/2$ and $\delta>0$ be arbitrary.
We can find $N(\delta)>0$ such that for any natural number $N\geq N(\delta)$ there exist a finite simplicial complex $K_N$ and
a $\delta$-embedding $\varphi_N: \left(Y, d^\prime_N\right) \to K_N$ satisfying 
\begin{equation}  \label{eq: dimension bound on K_N}
   \frac{\dim K_N}{N} < \mdim(Y, S) + \delta.
\end{equation}
For each natural number $n$ with $N_n \geq N(\delta)$ we consider a map
\[  \pi_n := \varphi_{N_n}\circ \pi \circ \psi_n: [0,1]^{M_n} \to K_{N_n}. \]
We use a famous result of topological dimension theory \cite[Theorem V.2]{Hurewicz--Wallman}:
The space $K_N$ topologically embeds in the Euclidean space $\mathbb{R}^{2\dim K_N +1}$.
(Since $K_N$ is a finite simplicial complex, this is a rather elementary fact.
See e.g. \cite[Lemma 2.6]{Gutman--Qiao--Tsukamoto}.)
Therefore we can consider that 
\[   K_{N_n} \subset \mathbb{R}^{2\dim K_{N_n} + 1}. \]
So $\pi_n$ is a continuous map from $[0,1]^{M_n}$ to $\mathbb{R}^{2\dim K_{N_n}+1}$.
We apply the waist inequality to this map. The rest of the argument is quite similar to the case (2).

By Corollary \ref{corollary: waist inequality for covering number} we can find $t\in K_{N_n}$ (depending on $n$) such that
\begin{equation}  \label{eq: waist inequality for the proof of main theorem again}
     \#\left(\pi_n^{-1}(t), \norm{\cdot}_\infty, \varepsilon\right) 
     \geq \frac{1}{8^{M_n}}\cdot \left(\frac{1}{\varepsilon}\right)^{M_n-2\dim K_{N_n}-1}. 
\end{equation}     
Pick $t^\prime \in Y$ (depending on $n$) with $\varphi_{N_n}(t^\prime) = t$.
We claim that 
\begin{equation} \label{eq: fiber of pi_n again}
   \psi_n\left(\pi_n^{-1}(t)\right) \subset \pi^{-1}\left(B_\delta\left(t^\prime, d^\prime_{N_n}\right)\right).
\end{equation}
Indeed, take $x\in \pi_n^{-1}(t)$. Then $\varphi_{N_n}\left(\pi\left(\psi_n(x)\right)\right) = \pi_n(x) = t = \varphi_{N_n}(t^\prime)$.
Since $\varphi_{N_n}$ is a $\delta$-embedding with respect to $d^\prime_{N_n}$, we have 
\[  d^\prime_{N_n}\left(\pi\left(\psi_n(x)\right), t^\prime\right) < \delta. \]
So $\pi\left(\psi_n(x)\right) \in B_\delta\left(t^\prime, d^\prime_{N_n}\right)$ and hence 
$\psi_n(x) \in \pi^{-1}\left(B_\delta\left(t^\prime, d^\prime_{N_n}\right)\right)$.
This has shown (\ref{eq: fiber of pi_n again}).

Since $\norm{x-y}_\infty \leq d_{N_n}\left(\psi_n(x), \psi_n(y)\right)$ $(x, y\in [0,1]^{M_n})$, 
\begin{align*}
   \#\left(\pi_n^{-1}(t), \norm{\cdot}_\infty, \varepsilon\right) & \leq \#\left(\psi_n\left(\pi_n^{-1}(t)\right), d_{N_n}, \varepsilon\right) \\
    &\leq \#\left(\pi^{-1}\left(B_\delta\left(t^\prime, d^\prime_{N_n}\right)\right), d_{N_n}, \varepsilon\right) \quad 
    \text{by (\ref{eq: fiber of pi_n again})}.
\end{align*}
By (\ref{eq: waist inequality for the proof of main theorem again}) we have 
\[  \sup_{y\in Y} \#\left(\pi^{-1}\left(B_\delta(y, d^\prime_{N_n})\right), d_{N_n}, \varepsilon\right) 
     \geq \frac{1}{8^{M_n}} \cdot \left(\frac{1}{\varepsilon}\right)^{M_n-2\dim K_{N_n} -1}. \]
Taking the logarithm and dividing it by $N_n \log(1/\varepsilon)$, we have
\begin{align*}
   \frac{\sup_{y\in Y} \log 
   \#\left(\pi^{-1}\left(B_\delta(y, d^\prime_{N_n})\right), d_{N_n}, \varepsilon\right)}{N_n \log(1/\varepsilon)} & \geq 
    \frac{M_n}{N_n} - \frac{2\dim K_n}{N_n} - \frac{1}{N_n} - \frac{M_n \log 8}{N_n \log (1/\varepsilon)} \\
    & \geq  \frac{M_n}{N_n} - 2\mdim(Y,S)-2\delta - \frac{1}{N_n} - \frac{M_n \log 8}{N_n \log (1/\varepsilon)}.
\end{align*}
Here we have used (\ref{eq: dimension bound on K_N}) in the second inequality.

We let $n\to \infty$. Noting $M_n/N_n\to s$, we get 
\[  \lim_{N\to \infty} \frac{\sup_{y\in Y} \log 
   \#\left(\pi^{-1}\left(B_\delta(y, d^\prime_{N})\right), d_{N}, \varepsilon\right)}{N \log(1/\varepsilon)} 
   \geq s- 2\mdim(Y, S) - 2\delta - \frac{s\log 8}{\log(1/\varepsilon)}.  \]
We let $\delta\to 0$ and then $\varepsilon\to 0$.
We conclude  $\lmdimm(\pi, T, d) \geq s-2\mdim(Y, S)$ by Lemma \ref{lemma: formula of conditional metric mean dimension}.
\end{proof}

\begin{remark}  \label{remark: proof of the main technical theorem} 
 There are two points in the above proof which might be able to be improved in a future:
 \begin{enumerate}
  \item In the above proof, we used the waist inequality in (\ref{eq: consequence of waist inequality in main technical theorem})
           and (\ref{eq: waist inequality for the proof of main theorem again}).
           We notice that they are slightly weaker than the original waist inequalities.
           In both the inequalities 
           (\ref{eq: consequence of waist inequality in main technical theorem}) 
           and (\ref{eq: waist inequality for the proof of main theorem again})
           the parameter $\varepsilon$ is fixed.
           On the other hand, it can vary in the original waist inequalities.
            It is an interesting problem to discover how to apply the full-power of the waist inequality 
            to the study of mean dimension.  A related question is the following: 
            For an equivariant continuous map $\pi:(X, T) \to (Y, S)$ between dynamical systems, 
            we define the \textbf{fiber-wise lower conditional metric mean dimension} by 
            \[ \lmdimm(\pi, T, d)_{\mathrm{fiber}} = 
            \sup_{y\in Y} \left\{\liminf_{\varepsilon\to 0} 
            \left(\limsup_{N\to \infty} \frac{\log \#\left(\pi^{-1}(y), d_N, \varepsilon\right)}{N\log (1/\varepsilon)}\right)\right\}. \]
            Can one replace $\lmdimm(\pi, T, d)$ with $\lmdimm(\pi, T, d)_{\mathrm{fiber}}$ in the statements (2) and (3) in 
            Theorem \ref{theorem: main technical theorem}?

  \item  The multiplicative factor $2$ of $2\mdim(Y, S)$ in the statement of Theorem \ref{theorem: main technical theorem} (3) comes form 
            the embedding $K_{N_n} \subset \mathbb{R}^{2\dim K_{N_n} + 1}$ used in the above proof.
            (See also Problem \ref{problem: reduce the factor two to one} in \S \ref{subsection: conditional metric mean dimension}.)
            If we would like to reduce the factor $2$ to a smaller value (hopefully, one), then we need to develop a method to bypass this embedding.
            It seems possible to find such a method if the local topology of the simplicial complex $K_{N_n}$ is not very complicated.
            For example, if $K_{N_n}$ is a smooth manifold then there exists a smooth map $f:K_{N_n}\to \mathbb{R}^{\dim K_{N_n}}$ 
            for which every fiber has cardinality at most $4\dim K_{N_n}$ \cite[p. 447]{Gromov10}. 
            Then we apply the waist inequality to $f\circ \pi_n$ and get the desired result.
            Another idea is to develop a waist inequality directly applicable to the map $\pi_n:[0,1]^{M_n}\to K_{N_n}$.
            But currently we do not know how to do it.
 \end{enumerate}  
\end{remark}

Theorems \ref{theorem: conditional metric mean dimension between two full shifts}, 
\ref{theorem: conditional metric mean dimension and mean dimension} and Proposition \ref{prop: analogue of metric Hurewicz}
in \S \ref{section: introduction} are immediate corollaries of the above main technical theorem.
They are all included in the following statement.

\begin{corollary}[$=$ Theorems \ref{theorem: conditional metric mean dimension between two full shifts}, 
\ref{theorem: conditional metric mean dimension and mean dimension} and Proposition \ref{prop: analogue of metric Hurewicz}]
Let $a$ be a natural number, and let $\left(\left([0,1]^a\right)^{\mathbb{Z}},\sigma\right)$ be the full-shift on the alphabet $[0,1]^a$.
We define a metric $D$ on it by 
\[   D\left((x_n)_{n\in \mathbb{Z}}, (y_n)_{n\in \mathbb{Z}}\right) = \sum_{n\in \mathbb{Z}} 2^{-|n|} \norm{x_n-y_n}_\infty. \]
  \begin{enumerate}
     \item  Let $b$ be a natural number. For any equivariant continuous map 
              $\pi:\left(\left([0,1]^a\right)^{\mathbb{Z}},\sigma\right) \to \left(\left([0,1]^b\right)^{\mathbb{Z}},\sigma\right)$ we have
              $\lmdimm\left(\pi, \sigma, D\right) \geq a-b$.
     \item Let $(Y, S)$ be a dynamical system. For any equivariant continuous map 
              $\pi:\left(\left([0,1]^a\right)^{\mathbb{Z}},\sigma\right) \to (Y, S)$ we have 
              $\lmdimm\left(\pi, \sigma, D\right) \geq a - 2\mdim(Y, S)$.
              
              Notice that this statement is logically contained in the next claim (3) (letting $Z$ to be one-point there). 
              But we separately state it for making the exposition easier to understand. 
              
      \item  Let $(Z,R)$ be a dynamical system with a metric $\rho$.
              We define a metric $d$ on the product $\left([0,1]^a\right)^{\mathbb{Z}}\times Z$ by 
              $d\left((x,z), (x^\prime,z^\prime)\right) = D(x,x^\prime) + \rho(z,z^\prime)$.
              Then for any dynamical system $(Y, S)$ and any equivariant continuous map 
              $\pi: \left(\left([0,1]^a\right)^{\mathbb{Z}} \times Z, \sigma\times R\right) \to (Y, S)$ we have
              $\lmdimm\left(\pi, \sigma\times R, d\right) \geq a - 2\mdim(Y, S)$.
  \end{enumerate}
\end{corollary}

\begin{proof}
For each natural number $n$ we define a map $\psi_n: \left([0,1]^a\right)^n \to \left([0,1]^a\right)^{\mathbb{Z}}$ by 
\[   \psi_n(x_0, x_1, \dots, x_{n-1})_k = \begin{cases} x_k & (0\leq k \leq n-1) \\ 0 & (\text{otherwise}) \end{cases}.     \]
This satisfies that for any $x, y\in  \left([0,1]^a\right)^n$
\[  \norm{x-y}_\infty \leq  D_n\left(\psi_n(x), \psi_n(y)\right).   \]
Then the statements (1) and (2) immediately follow from Theorem \ref{theorem: main technical theorem} with the parameters
$N_n := n$ and $M_n := an$.

For (3), we fix a point $p\in Z$ and consider a map 
\[  \psi_n^\prime:\left([0,1]^a\right)^n \to \left([0,1]^a\right)^{\mathbb{Z}} \times Z, \quad 
    x\mapsto \left(\psi_n(x), p\right).     \]
This satisfies that for any $x, y\in  \left([0,1]^a\right)^n$
\[  \norm{x-y}_\infty \leq  d_n\left(\psi^\prime_n(x), \psi^\prime_n(y)\right).     \]
We use Theorem \ref{theorem: main technical theorem} with these maps $\psi^\prime_n$ and get the statement (3).
\end{proof}

\section{Construction of a minimal closed subset}   \label{section: construction of a minimal subset}

We prove Theorem \ref{theorem: conditional metric mean dimension and minimal system} in this section.
The construction of a minimal closed subset itself is the same as that of \cite[Proposition 3.5]{Lindenstrauss--Weiss}.
It uses an idea of \textit{block-system} reviewed below.

\subsection{Review of block-system}  \label{subsection: review of block-system}

Let $a$ be a natural number. 
As before, we define a metric $D$ on the full-shift $\left([0,1]^a\right)^{\mathbb{Z}}$ by 
\[   D\left((x_n)_{n\in \mathbb{Z}}, (y_n)_{n\in \mathbb{Z}}\right) = \sum_{n\in \mathbb{Z}} 2^{-|n|} \norm{x_n-y_n}_\infty. \]
Let $x$ be a point in $\left([0,1]^a\right)^{\mathbb{Z}}$.
For integers $\ell<m$ we set 
\[ x|_\ell^m = (x_\ell, x_{\ell+1}, \dots, x_{m-1}, x_m)  \in \left([0,1]^a\right)^{m-\ell+1}. \]

Let $N$ be a natural number, and
let $K$ be a closed subset of $\left([0,1]^a\right)^N$.
We define a \textbf{block-system} $X(K)$ in $\left([0,1]^a\right)^{\mathbb{Z}}$ by 
\[  X(K) = \left\{x\in \left([0,1]^a\right)^{\mathbb{Z}} \middle|\, \exists \ell \in \mathbb{Z}: \forall n\in \mathbb{Z}:
                x|_{\ell+nN}^{\ell+(n+1)N-1} \in K\right\}. \]
This is a shift-invariant closed subset of $\left([0,1]^a\right)^{\mathbb{Z}}$.

\begin{lemma}   \label{lemma: metric mean dimension of block-system}
\[  \umdimm\left(X(K),\sigma, D\right) \leq \frac{\udimm(K, \norm{\cdot}_\infty)}{N}. \]
Here $\sigma$ is the shift-map, and $\udimm(K, \norm{\cdot}_\infty)$ is the upper Minkowski dimension of $K$ with respect to 
the $\ell^\infty$-distance.
\end{lemma}

\begin{proof}
For $0\leq \ell \leq N-1$ we set 
\[  X_\ell = \left\{x\in \left([0,1]^a\right)^{\mathbb{Z}}\middle|\, \forall n\in \mathbb{Z}: x|_{\ell+nN}^{\ell+(n+1)N-1} \in K\right\}. \]
We have 
\[  X(K) = X_0 \cup X_1 \cup X_2 \cup \dots \cup X_{N-1}. \]
Let $\varepsilon$ be a positive number. Take $m = m(\varepsilon)>0$ satisfying 
$\sum_{|n|>m} 2^{-|n|} < \varepsilon/2$.

Let $L$ be an arbitrary natural number. We set 
\[  s = \left\lfloor \frac{-\ell-m}{N}\right\rfloor, \quad 
     t = \left\lceil \frac{L-\ell+m}{N}\right\rceil -1. \]
We define $f_\ell: X_\ell \to K^{t-s+1}$ by 
\[  f_\ell(x) = \left(x|_{\ell+nN}^{\ell+(n+1)N-1}\right)_{s\leq n \leq t}. \]
Notice that 
\[  \bigcup_{n=s}^t [\ell+nN, \ell+(n+1)N) = [\ell+sN, \ell+(t+1)N) \supset [-m, L+m). \]
Hence, for any $x,y\in X_\ell$
\[   D_L\left(x, y\right) < 3 \norm{f_\ell(x)-f_\ell(y)}_\infty + \frac{\varepsilon}{2}. \]
Then 
\[  \#\left(X_\ell, D_L, \varepsilon\right) \leq \#\left(K^{t-s+1}, \norm{\cdot}_\infty, \frac{\varepsilon}{9}\right)   
     \leq \left(\#\left(K, \norm{\cdot}_\infty, \frac{\varepsilon}{9}\right)\right)^{t-s+1}. \]
We have 
\[  t-s+1 \leq \frac{L-\ell+m}{N} - \frac{-\ell-m}{N}+1+1 \leq \frac{L+2m}{N}+2. \]
Hence 
\[  \#\left(X_\ell, D_L, \varepsilon\right) \leq \left(\#\left(K, \norm{\cdot}_\infty, \frac{\varepsilon}{9}\right)\right)^{\frac{L+2m}{N}+2}. \]
Therefore 
\[  \#\left(X(K), D_L, \varepsilon\right) \leq N\cdot \left(\#\left(K, \norm{\cdot}_\infty, \frac{\varepsilon}{9}\right)\right)^{\frac{L+2m}{N}+2}. \]
Taking the logarithm and dividing it by $L\log(1/\varepsilon)$, we get
\[  \frac{\log \#\left(X(K), D_L, \varepsilon\right)}{L\log(1/\varepsilon)}
     \leq \frac{\log N}{L\log(1/\varepsilon)} + \left(\frac{1}{N}+\frac{2m}{NL}+\frac{2}{L}\right) 
            \frac{\log \#\left(K, \norm{\cdot}_\infty, \frac{\varepsilon}{9}\right)}{\log (1/\varepsilon)}. \]
We let $L\to +\infty$ and then let $\varepsilon\to 0$.
We conclude that 
\[  \umdimm\left(X(K),\sigma, D\right) \leq \frac{\udimm(K, \norm{\cdot}_\infty)}{N}. \]
\end{proof}

\subsection{Proof of Theorem \ref{theorem: conditional metric mean dimension and minimal system}}
\label{Proof of theorem of conditional metric mean dimension and minimal system}

Now we prove Theorem \ref{theorem: conditional metric mean dimension and minimal system}.
We write the statement again.

\begin{theorem}[$=$ Theorem \ref{theorem: conditional metric mean dimension and minimal system}]
Let $a$ be a natural number and $s$ a real number with $0\leq s < a$.
There exists a minimal closed subset $X$ of $\left(\left([0,1]^a\right)^{\mathbb{Z}}, \sigma\right)$ satisfying the following 
two conditions.
   \begin{enumerate}
     \item $\mdim\left(X, \sigma\right) = \mdimm\left(X, \sigma, D\right) = s$.
     \item For any natural number $b$ and any equivariant continuous map 
              $\pi:(X,\sigma) \to \left(\left([0,1]^b\right)^{\mathbb{Z}},\sigma\right)$ we have 
              $\lmdimm\left(\pi, \sigma, D\right) \geq s-b$.
   \end{enumerate}
\end{theorem}

\begin{proof}
As we already noted, the construction of $X$ is the same as that of \cite[Proposition 3.5]{Lindenstrauss--Weiss}.
The difference is that we analyze the resulting minimal closed subset by using the waist inequality.

We take and fix a sequence of rational numbers $0<r_n<1$ $(n\geq 1)$ satisfying 
\[  a \prod_{n=1}^\infty (1-r_n) = s. \]
Set $N_1 = 1$ and $K_1 = [0,1]^a$.
Starting from these, we inductively construct an increasing sequence of natural numbers $N_n$ and closed subsets
$K_n\subset \left([0,1]^a\right)^{N_n}$ $(n\geq 1)$ satisfying the following conditions.
\begin{enumerate}
   \item[(i)] There are natural numbers $p_n > q_n$ $(n\geq 1)$ larger than one satisfying $N_{n+1} = p_n N_n$ and $r_n = q_n/p_n$.
   \item[(ii)] $K_{n+1}\subset \left(K_n\right)^{p_n}$. Moreover,
                  there is a point $(w^{(n)}_1, w^{(n)}_2, \dots, w^{(n)}_{q_n})$ in $\left(K_n\right)^{q_n}$ such that 
   \begin{itemize}
     \item $K_{n+1} = (K_n)^{p_n-q_n} \times \{(w^{(n)}_1, w^{(n)}_2,\dots, w^{(n)}_{q_n})\}$.
     \item For any $x\in K_n\times K_n$ there exists $1\leq k < q_n$ satisfying 
              \[   \norm{x-(w^{(n)}_k, w^{(n)}_{k+1})}_\infty < \frac{1}{n}.  \]   
   \end{itemize}
\end{enumerate}
A moment thought shows that such an inductive construction works well.

We consider the block-system $X_n = X(K_n)$ associated with $K_n$. We set 
\[  X  =  \bigcap_{n=1}^\infty X_n \subset \left([0,1]^a\right)^{\mathbb{Z}}. \]
By the above condition (ii) the set $X$ is minimal.
We will show that it satisfies the required conditions (1) and (2).

We define a sequence of natural numbers $M_n$ $(n\geq 1)$ by 
\[  M_1 = 1, \quad M_{n+1} = (p_n-q_n) M_n \quad  (n\geq 1). \]
We have 
\[  M_{n+1} = p_n\left(1-\frac{q_n}{p_n}\right) M_n = p_n (1-r_n) M_n = p_n N_n (1-r_n) \frac{M_n}{N_n} = N_{n+1} (1-r_n)\frac{M_n}{N_n}. \]
So $\frac{M_{n+1}}{N_{n+1}} = (1-r_n) \frac{M_n}{N_n}$ and hence 
\[  \frac{M_{n+1}}{N_{n+1}} = \prod_{k=1}^n (1-r_k) \to \frac{s}{a} \quad (n\to \infty). \]

From the construction $K_n$ is naturally homeomorphic to the cube $([0,1]^a)^{M_n}$.
More precisely, $(K_n, \norm{\cdot}_\infty)$ is isometric to $\left(([0,1]^a)^{M_n}, \norm{\cdot}_\infty \right)$.
In particular, 
\[  \dim_{\mathrm{M}}\left(K_n, \norm{\cdot}_\infty\right) = a M_n. \]
By Lemma \ref{lemma: metric mean dimension of block-system}
\[  \umdimm\left(X_n, \sigma, D\right) \leq \frac{\dim_{\mathrm{M}}\left(K_n, \norm{\cdot}_\infty\right)}{N_n}  
     = \frac{a M_n}{N_n}. \]
Thus 
\[  \umdimm\left(X, \sigma, D\right) \leq  \lim_{n\to \infty} \frac{a M_n}{N_n} = s.  \]
So we have $\umdimm\left(X, \sigma, D\right) \leq s$.

It also follows from the construction that there exists a continuous map 
$\psi_n: \left([0,1]^a\right)^{M_n} \to X$ satisfying 
\[   \norm{\psi_n(x)|_0^{N_n-1} - \psi_n(y)|_0^{N_n-1}}_\infty = \norm{x-y}_\infty \quad 
      (x, y \in \left([0,1]^a\right)^{M_n}). \]
In particular 
\[   \norm{x-y}_\infty \leq  D_{N_n}\left(\psi_n(x), \psi_n(y)\right). \]
Then we can apply to the maps $\psi_n$ Theorem \ref{theorem: main technical theorem} with the parameters $a M_n$ and $N_n$.
Noting $\lim_{n\to \infty} \frac{a M_n}{N_n} = s$, we get the following conclusions:
\begin{itemize}
  \item $\mdim(X, \sigma) \geq  s$.
  \item For any natural number $b$ and 
          any equivariant continuous map $\pi: (X, \sigma) \to \left(\left([0,1]^b\right)^{\mathbb{Z}}, \sigma\right)$ we have 
          $\lmdimm\left(\pi, \sigma, D\right) \geq s - b$.
\end{itemize}
Since mean dimension is bounded from above by metric mean dimension \cite[Theorem 4.2]{Lindenstrauss--Weiss}, 
we have 
\[  s\leq \mdim(X, \sigma)  \leq  \lmdimm\left(X, \sigma, D\right)  \leq  \umdimm\left(X, \sigma, D\right) \leq  s. \]
Thus $\mdim(X, \sigma) = \mdimm(X, \sigma, D) = s$.
So $(X, \sigma, D)$ satisfies all the required conditions.
\end{proof}

\end{document}